\documentclass[12pt,a4paper]{amsart}
\usepackage{color}
\usepackage{amssymb}

\def\Q{\mathbb{Q}}
\def\Z{\mathbb{Z}}

\def\C{\mathbb{C}}

\def\<{\langle}
\def\>{\rangle}

\newtheorem{thmA}{Theorem}

\newtheorem{theorem}{Theorem}
\newtheorem{lemma}[theorem]{Lemma}

\newenvironment{pf}{\par\medskip\noindent\textit{Proof.}~}{\hfill $\square$\par\medskip}

\begin{document}

\catcode`\@=11
\def\serieslogo@{\relax}
\def\@setcopyright{\relax}
\catcode`\@=12

\def\jump{\vskip 0.3cm}
\title[A remark on double cosets]{A remark on double cosets}

\author[Howie]{James Howie }
\address{ James Howie\\
Department of Mathematics and Maxwell Institute for Mathematical Sciences\\
Heriot--Watt University\\
Edinburgh EH14 4AS }
\email{ jim@ma.hw.ac.uk}

\begin{abstract}
If a soluble group $G$ contains two finitely generated abelian
subgroups $A,B$ such that the number of double cosets $AgB$ is finite,
then $G$ is shown to be virtually polycyclic.
\end{abstract}

\maketitle

\section{Introduction}

Suppose $G$ is a group, and $A,B$ are subgroups of $G$ that are, in some sense,
`small'.  If the set $A\backslash G/B$ of double cosets $AgB$ ($g\in G$) of $A$ and $B$ in $G$ is finite,
does it necessarily follow that $G$ is also `small'?

Questions of this nature have been the subject of some previous investigations in
the case where $A=B$ \cite{Niblo,J1,J2}.  Under a variety of suitable conditions, it
turns out that the hypotheses imply that $A$ has finite index in $G$.  This is
related to questions of Dunwoody \cite{Dun} and of Neumann and Rowley \cite{NR}.

The general question arises in connection with the study of homological finiteness
properties of subdirect products of limit groups \cite{bh2,bhms2}.  For example
in \cite{bh2} it is shown that a certain rational homology group $H_n(G,\Q)$
contains a
$\Q$-linearly independent set in one-to-one correspondence with a set
$A\backslash G/B$ of double cosets, where $A$ is a cyclic subgroup and $B$ a
finitely generated abelian group.  In the context of \cite{bh2}, the group $G$
is readily shown not to be `small' enough for such a set to be finite, from which
one deduces that it is not of type $FP_n(\Q)$.

The purpose of this short note is to prove the following.

\begin{thmA} \label{theoremC}
Let $G$ be a virtually soluble group, and let
$A,B<G$ be finitely generated virtually abelian subgroups of $G$.
If there are only finitely many double cosets $AgB$ ($g\in G$), then
$G$ is virtually polycyclic.
\end{thmA}

I first proved this result in the case where $G$ is (virtually) nilpotent-by-abelian, and
$A$ is cyclic, with the intention of using it in the proof of the main theorem of
\cite{bhms2} (that a subdirect product of $n$ limit groups has type $FP_n(\Q)$
only if it is virtually a direct prodcut of $n$ or fewer limit groups).
In the event, an alternative more direct approach was developed for \cite{bhms2}
which avoided the use of double cosets. Nevertheless the result seems of some
interest in relation to the various references cited above.

The special case $A=B$ is connected to the question of orbits of groups
of automorphisms: if $G=N\rtimes A$, then $A\backslash G/A$ is in one-to-one
correspondence with the set of orbits of the action of $A$ on $N$.  A partial
analogue of this correspondence is a useful tool for the study of the more
general case.  Suppose that $G=N\rtimes A=N\rtimes B$ -- in particular $A\cong B$
via a preferred isomorphism $\theta:A\to B$ such that $a\theta(a)^{-1}\in N$ for all
$a\in A$.  Restricting to the situation where $N$ is abelian, the map
$\delta:A\to N$, $\delta(a):=a\theta(a)^{-1}$ is a derivation, so gives rise
to a $\Z A$-module homomorphism from $I$ to $N$ via the rule $1-a\mapsto\delta(a)$
(where $I$ denotes the augmentation ideal in the group ring $\Z A$).  We can
also associate to this set-up an action of $A$ on $N$ by {\em affine transformations}
$$a*\nu:=a\nu\theta(a)^{-1}=a\cdot\nu+\delta(a),$$
where $a\cdot\nu=a\nu a^{-1}$ is the standard conjugation action of $A$ on $N$.
It is not difficult to check that the set $A\backslash G/B$ is in one-to-one
correspondence with the set of orbits of this affine action of $A$ on $N$.

Thus a special case of Theorem \ref{theoremC} says that, if a finitely generated virtually
abelian group $A$ acts affinely on an abelian group $N$ with finitely many orbits, then
$N$ is finitely generated.

\begin{thmA}\label{theoremA}
Let $A$ be a finitely generated virtually abelian group, $N$ a $\Z A$-module, and
$\delta:A\to N$ a derivation.  If the resulting affine action of $A$ on $N$ has only
finitely many orbits, then $N$ is finitely generated as an abelian group.
\end{thmA}

This should be compared with a theorem of Jabara \cite{J2}, who proves that if
the action of $A$ by {\em automorphisms} on a (not necessarily abelian) group $N$
has only finitely many orbits, then $N$ is finite.

Theorem \ref{theoremA} is a key step in our proof of Theorem \ref{theoremC}.
We prove this in Section \ref{aff} below, and then complete the proof of Theorem \ref{theoremC}
in Section \ref{doublecosets}.

I am grateful to my co-authors of \cite{bhms2}, Martin Bridson, Chuck Miller
and Hamish Short, for a fruitful and enjoyable collaboration, and for their encouragement
to publish these double-coset results separately.
I was greatly encouraged while working on this problem by the interest in it shown
by Karl Gruenberg, in the course of a number of conversations and letters (which
culminated in his independent solution of a special case).  Karl was a great source of
inspiration to me, as he was to many mathematicians.  This paper is dedicated to his
memory.

\section{Affine actions}\label{aff}

Let $A$ be a group, and let $N$ be a (left) $\Z A$-module.
Recall that a map $\delta:\to N$ is a {\em derivation} if it satisfies the
Leibnitz rule $$\delta(a_1a_2)=\delta(a_1)+a_1\cdot\delta(a_2)$$
for all $a_1,a_2\in A$.  Here we use $a\cdot\nu$, $a\in A$, $\nu\in N$,
to denote the module action of $A$ on $N$.

Given a derivation $\delta:A\to N$, we define the corresponding
{\em affine action} of $A$ on $N$ by
$$a*\nu:=a\cdot\nu+\delta(a).$$

Note that $$a_1*(a_2*\nu)=a_1\cdot(a_2\cdot\nu+\delta(a_2))+\delta(a_1)$$
$$=a_1\cdot(a_2\cdot\nu)+(\delta(a_1)+a_1\cdot\delta(a_2))=(a_1a_2)\cdot\nu+\delta(a_1a_2),$$
so this is a genuine action of $A$ on the underlying set of $N$, although it does
not in general respect the group structure of $N$.  In the special case where $\delta$
is the zero derivation, the affine action is just the standard action of $A$ by group
automorphisms on $N$.

It is well-known (see for example \cite[IV.2, Exercise 2]{ksbrown})
that the abelian group $Der(A,N)$ of derivations
$A\to N$ (where the binary operation is pointwise addition in $N$)
is isomorphic to $Hom_{\Z A}(I,N)$, where $I$ is the augmentation ideal
in $\Z A$.  Here $\phi:I\to N$ corresponds to $\delta:A\to N$
defined by $\delta(a)=\phi(1-a)$.  If $J$ is any left ideal in $\Z A$,
there is a canonical $\Z A$-homomorphism $I\to \Z A/J$ given by
$x\mapsto x+J$.  We call the corresponding derivation
$a\mapsto 1-a+J$ and affine action $a*(x+J):=(ax+1-a)+J$
the {\em canonical derivation} and {\em canonincal affine action}
respectively.

In all the above constructions, one can replace $\Z$ by another commutative
ring.  In practice, we will require only the finite fields $\Z_p$ and
the field $\Q$ of rationals for this purpose.

\begin{lemma}\label{affine_poly}
Let $R$ be either $\Z$ or $\Z_p$ for some prime $p$, let $A$ be a finitely generated
abelian group, $I$ the augmentation ideal in $RA$, and $J\subset I$ another ideal in
$RA$.  Suppose that $I/J$ is the union of finitely many $A$-orbits under the
canonical affine action of $A$ defined by
$$a*(x+J):=(ax+1-a)+J,~~(a\in A,~x\in I).$$
Then for each $a\in A$ there exists a nonzero polynomial $f(X)\in R[X]$
such that $f(a)\in J$.
\end{lemma}

\begin{pf}
Let $K$ be an algebraically closed field containing $R$.
For each prime $q$, let $f_q(t)=(t-1)^q$.  In particular, $f_q(a)\in I$ for
all $a\in A$.  Let $\Pi$ denote the infinite set of primes (if $R=\Z$), or
the set of primes $\ne p$ (if $R=\Z_p$).  Since $I/J$ is the union of finitely
many orbits under the affine action, there are two distinct primes $q_1,r_1\in\Pi$
and an element $a_1\in A$ such that $f_{q_1}(a)-a_1*f_{r_1}(a)\in J$.
Applying the same argument to the infinite set $\Pi\smallsetminus\{q_1,r_1\}$, there
are two more primes $q_2,r_2$ and $a_2\in A$ such that $f_{q_2}(a)-a_2*f_{r_2}(a)\in J$.
Iterating this argument yields an infinite sequence $q_1,r_1,q_2,r_2,q_3,r_3,\dots$
of pairwise distinct primes ($\ne p$ in the case $R=\Z_p$) and a sequence
$a_1,a_2,\dots$ of elements of $A$, such that $f_{q_j}(a)-a_j*f_{r_j}(a)\in J$ for
each $j$.

Since $A$ is finitely generated, there is a natural number $n$ and a
nonzero $n$-tuple $(t_1,\dots,t_n)$
of integers such that
\begin{equation}\label{lindep}
a_1^{t_1}\cdot a_2^{t_2} \cdots a_n^{t_n}=1~~\mathrm{in}~A.
\end{equation}

Note that $f_{r_j}(a)-a_j^{-1}*f_{q_j}(a)\in J$ for each $j$, so by interchanging
$q_j$ and $r_j$ for each $j$ in  suitable subset of $\{1,\dots,n\}$, we can
arrange that each $t_j$ is non-negative in (\ref{lindep}).

Now the equations
$$(a-1)^{q_j} - 1 = f_{q_j}(a)-1 = a_j*f_{r_j}(a)-1$$
$$= a_j\left( f_{r_j}(a)-1\right) =a_j\left( (a-1)^{r_j}-1\right)$$
in $I/J$ for $j=1,\dots,n$ give rise to an equation
$$\prod_{j=1}^n \left( (a-1)^{q_j} - 1\right)^{t_j} = \prod_{j=1}^n a_j^{t_j} \prod_{j=1}^n \left( (a-1)^{r_j} - 1\right)^{t_j} = \prod_{j=1}^n \left( (a-1)^{r_j} - 1\right)^{t_j}$$
in $I/J$.  If $\alpha\in K$ is a
primitive $q_1$-th root of unity, then
$\alpha+1$ is a root of $g_1(X):=\prod_{j=1}^n \left( (X-1)^{q_j} - 1\right)^{t_j}$
but not
of $g_2(X):=\prod_{j=1}^n \left( (X-1)^{r_j} - 1\right)^{t_j}$.  Hence $f(X):=g_1(X)-g_2(X)\ne 0$ in $R[X]$,
while $f(a)\in J$ as required.
\end{pf}

\begin{lemma}\label{restrict}
Let $A$ be a group, $M$ a $\Z A$-module, $M_0$ a submodule of $M$,
and $\delta:A\to M$ a derivation.  Then
\begin{enumerate}
 \item $A_0:=\delta^{-1}(M_0)$ is a subgroup of $A$;
\item the affine action $a*x=ax+\delta(a)$ of $A$ on $M$ restricts
to an affine action of $A_0$ on $M_0$, and induces an affine action of
$A$ on the quotient module $M/M_0$;
\item the orbits of the affine action of $A_0$ on $M_0$ are
the nonempty intersections of $M_0$ with the orbits of the affine
action of $A$ on $M$; and
\item the orbits of the affine action of $A$ on $M/M_0$ are the images
under the natural epimorphism $M\to M/M_0$ of the
orbits of the affine action of $A$ on $M$.
\end{enumerate}

\end{lemma}

\begin{pf}
\begin{enumerate}
 \item It follows immediately from the Leibnitz rule for derivations
that $A_0$ is a subgroup of $A$.
\item By definition, if $a\in A_0$ and $x\in M_0$ then $a*x\in M_0$,
so the affine action of $A$ on $M$ restricts to an affine action of
$A_0$ on $M_0$.
If $a\in A$, $x\in M$ and $y\in M_0$, then $a*(x+y)-a*(x)=a\cdot y\in M_0$,
so there is a well-defined induced action of $A$ on $M/M_0$.
\item Two elements $x,y\in M_0$ are in the same orbit of the
affine $A$-action on $M$ if and only if $y=a*x=a\cdot x+\delta(a)$ for
some $a\in A$.  But $\delta(a)=y-a\cdot x\in M_0$ means that $a\in A_0$.
\item Let $x,y\in M$ and $a\in A$.  If $y+M_0=a*(x+M_0)=(a*x)+M_0$
then $\exists~y'\in y+M_0$ with $a*x=y'$.  It follows that the orbits of the
affine $A$-action on $M/M_0$ are precisely the images of the orbits of the
affine $A$-action on $M$.
\end{enumerate}

\end{pf}

The next few lemmas take care of special cases of Theorem \ref{theoremA}.

\begin{lemma}\label{tors}
Let $A$ be a finitely generated virtually abelian group, $N$ a $\Z A$-module, and
$\delta:A\to N$ a derivation.  If the resulting affine action of $A$ on $N$ has only
finitely many orbits, then the $\Z$-torsion subgroup $N_0$ of $N$ is finite.
\end{lemma}

\begin{pf}
Passing to a subgroup of finite index if necessary, we may assume that
$A$ is free abelian of finite rank.
Now $N_0$ is a characteristic subgroup of $N$,
so is a $\Z A$-submodule.  By Lemma \ref{restrict}, $N_0$ is the union of finitely
many $A_0$ orbits under the restricted affine action of $A_0=\delta^{-1}(N_0)$
on $N_0$.  In particular, $N_0$ is finitely
generated as a $\Z A_0$-module.  Since it is also a $\Z$-torsion module, it
has finite exponent.  Arguing by induction on the number of prime factors of this
exponent, we will assume that $pN_0$ is finite for some prime $p$, and show that
$N_1=N_0/pN_0$ is also finite.

But $N_1$ is a finitely generated $\Z_p A_0$ module, which is the union
of finitely many $A_0$-orbits under the induced affine action.
If $I_0$ is the augmentation ideal in $\Z_pA_0$ and $\phi_0:I_0\to N_1$
the homomorphism corresponding to the derivation $\delta:A_0\to N_1$,
then the kernel of $\phi_0$ is an ideal $J_0$ of $\Z_pA_0$ with $J_0\subset I_0$.

Let $\{a_1,\dots,a_k\}$ be a basis for $A_0$ as a free abelian group.
Then by Lemma \ref{affine_poly} there are nonzero polynomials
$f_1(X),\dots,f_k(X)\in\Z_p[X]$ such that $f_1(a_1),\dots,f_k(a_k)\in J_0$.
It follows that $I_0/J_0$ is finite-dimensional as a $\Z_p$-vector space, and
hence finite.  Thus $\phi_0(I_0)$ is finite.

On the other hand, $N_1/\phi_0(I_0)$ is the union of finitely many orbits
under the induced affine action of $A_0$, by Lemma \ref{restrict}.
But this action is just the action by scalar multiplication, since
$\delta(A_0)\subset\phi_0(I_0)$.  It then follows from \cite[Theorem 1]{J2}
that $N_1/\phi_0(I_0)$ -- and hence also $N_1$ -- is finite, as required.
\end{pf}

\begin{lemma}\label{mult}
Let $A$ be a finitely generated virtually abelian group, $N$ a $\Z A$-module, and
$\delta:A\to N$ a derivation.  If the resulting affine action of $A$ on $N$ has only
finitely many orbits, then the quotient of $N$ by the submodule generated by $\delta(A)$ is finite.
\end{lemma}

\begin{pf}
Let $N_0$ denote the $\Z A$-submodule generated by $\delta(A)$.
Then the induced affine action of $A$ on $N/N_0$ is just the action by scalar
multiplication (since $\delta(A)\subset N_0$).  By Lemma \ref{restrict},
$N/N_0$ is the union of finitely many $A$-orbits under this action,
so by \cite[Theorem 1]{J2} it follows that $N/N_0$ is finite.
\end{pf}

\begin{lemma}\label{trans}
Let $A$ be a finitely generated virtually abelian group, $N$ a $\Z A$-module, and
$\delta:A\to N$ a derivation.  If the resulting affine action of $A$ on $N$ has only
finitely many orbits, then $N/IN$ is finitely generated as a $\Z$-module, where $I$
is the augmentation ideal in $\Z A$.
\end{lemma}

\begin{pf}
Since the module $A$-action on $N/IN$ is trivial, the resulting $A$-action is by translations:
$a*x=a\cdot x+\delta(a)=x+\delta(a)$.  The hypotheses then mean that $\delta:A\to N/IN$ is
a group homomorphism whose image has finite index in the additive group $N/IN$.  The result follows
from the fact that $A$ is finitely generated.
\end{pf}

\begin{theorem}[= Theorem \ref{theoremA}]
Let $A$ be a finitely generated virtually abelian group, $N$ a $\Z A$-module, and
$\delta:A\to N$ a derivation.  If the resulting affine action of $A$ on $N$ has only
finitely many orbits, then $N$ is finitely generated as a $\Z$-module.
\end{theorem}

\begin{pf}
Clearly, we may replace $A$ by a finite-index subgroup without changing
the hypotheses or the conclusion of the theorem.  Hence we may assume without loss of
generality that $A$ is free abelian of finite rank.

Now let $N_0$ be the $\Z$-torsion subgroup of $N$.  Then $N_0$ is finite, by
Lemma \ref{tors}.  Lemma \ref{restrict} thus reduces us to the case where
$N$ is $\Z$-torsion-free.

By Lemma \ref{mult}, the submodule of $N$ generated by $\delta(A)$
has finite index in $N$, so there is no loss of generality in assuming
that $N$ is generated by $\delta(A)$.

Let $\phi:I\to N$ denote the $\Z A$-module homomorphism corresponding to $\delta$.
Then $N\cong I/J$ for some ideal $J$ of $\Z A$ with $J\subset I$.

By Lemma \ref{trans}, the quotient $N/IN\cong I/(I^2+J)$ is finitely generated
as a $\Z$-module.  Moreover, by Lemma \ref{restrict}, $IN$ is the union of finitely
many $A_0$-orbits under the restriction of the affine action to 
$A_0=\delta^{-1}(IN)$, the kernel of the group homomorphism $A\to N/IN$.
If $A_0$ has strictly smaller rank than $A$, then the inductive hypothesis
asserts that $IN$ is finitely generated as a $\Z$-module, and the result follows.
Thus we are reduced to the case where, $A/A_0$ is finite -- and hence $N/IN$
is finite.

If $\{b_1,\dots,b_t\}$ is a basis for $A$ as an abelian group,
then Lemma \ref{affine_poly} says there are nonzero polynomials $g_1,\dots,g_t\in\Z[X]$
such that $g_i(b_i)\in J$ for $i=1,\dots,t$.  It follows that $\Q\otimes_\Z N$ has
finite dimension as a $\Q$-vector space.

Now $\Q\otimes_\Z N=\Q I/\Q J$, where $\Q I=\Q\otimes_\Z I$ is the augmentaion ideal
in $\Q A$, and $\Q J$ is the ideal of $\Q A$ generated by $J$.  Clearly
$\mathfrak{m}:=\Q I/\Q J$ is a maximal ideal in the ring $\Lambda:=\Q A/\Q J$.
Moreover, $\mathfrak{m}/\mathfrak{m}^2=0$, since $I/(I^2+J)=N/IN$ is finite.

If $\Lambda$ is a local ring, then $\mathfrak{m}$ is the unique maximal ideal
in $\Lambda$, and $\mathfrak{m}/\mathfrak{m}^2=0$, so Nakayama's Lemma (see for
example \cite[\S 2.8]{Reid}) implies
that $\mathfrak{m}=0$.  But $\mathfrak{m}=\Q\otimes_\Z N$ with $N$ torsion-free
as a $\Z$-module, so it follows that $N=0$.

Hence we may assume that $\Lambda$ is not local.  In other words, there is a maximal
ideal other than $\mathfrak{m}$ in $\Lambda$.  Thus there is a $\Q$-algebra
epimorphism $\psi$ from $\Lambda$ onto a field $K$, such that $\psi(I)\ne\{0\}$.
Note that $\Q A/\Q I\cong\Q$ and $\Q I/\Q J \cong \Q\otimes_\Z N$
are both finite-dimensional over $\Q$, so $\Lambda\cong \Q A/\Q J$ is finite-dimensional
over $\Q$.
It follows that $K$ is finite-dimensional over $\Q$, and so $K$ is an algebraic number field.

Let $d$ be the degree of $K$ over $\Q$.  Since $K$ is spanned by
$\psi(\Z A)=\Z+\psi(I)$ as a vector space over $\Q$,
we can choose $x\in I$ such that $K=\Q[\psi(x)]$.
Indeed, multiplying $x$ by a suitable positive integer,
we may assume in addition that $\psi(x)$ is an algebraic integer in $K$.

By hypothesis, there exists a finite subset $U\subset I$ with the
property that $I=\{a*u+y,a\in A,u\in U, y\in J\}$, where $a*u$
denotes the canonical affine action $a*u:=au+1-a$.  Note that
$a*u-1=a(u-1)$.  

Now define $\nu:I\to\Q$ by $\nu(z)=N_{K|\Q}(\psi(z)-1)$
for all $z\in I$, where $N_{K|\Q}$ denotes the norm of the algebraic extension $K|\Q$.

If $\xi_1,\dots,\xi_d$ denote the conjugates of $\psi(x)$
in $\C$, then $\xi_j\ne 0$ for each $j$.  Moreover,
for any positive integer $n$, the conjugates of $\psi(nx)$
are $n\xi_1,\dots,n\xi_d$.  Thus, replacing $x$ by a sufficiently large integer multiple, we may assume
that $|\xi_j|>2$ for each $j$.

Hence we have, for each integer $n\ge 1$,
$$\nu(nx)=N_{K|\Q}(\psi(nx)-1)=(n\xi_1-1)\cdots (n\xi_d-1),$$
an integer of abslolute value greater than $1$.
This expression also shows that
$$\nu(nx)=(-n)^d \mu(1/n),$$
where $$\mu(X)=(X-\xi_1)\cdots(X-\xi_d)\in\Z[X]$$ is the minimal polynomial of the algebraic integer $\psi(x)$.  In particular, $\nu(nx)$ is
congruent to $(-1)^d$ modulo $n$, so is coprime to $n$.
Since $\nu(nx)\ne\pm 1$, it has at least one prime divisor, and none of its prime divisors divide $n$.

Since we can do this for arbitrary integers $n$,
we deduce that there are infinitely many primes $p$
with the property $(\exists~z\in I)$,
$\nu(z)$ is a nonzero integer divisible by $p$.

On the other hand, every element of $I$ has the form
$a*u+y$ for some $a\in A$, $u\in U$ and $y\in J$,
and $$\nu(a*u+y)=N_{K|\Q}(\psi(au-a)=N_{K|\Q}(\psi(a)(\psi(u)-1))$$
$$=N_{K|\Q}(\psi(a))N_{K|\Q}(\psi(u)-1)=N_{K|\Q}(\psi(a))\nu(u).$$
It follows that, if $\nu(a*u+y)$ is a nonzero integer, then its prime divisors belong to a finite set -- namely
those primes that divide the numerator or denominator
of $N_{K|\Q}(\psi(a))$ for some $a$ in a fixed basis of $A$, together with
those primes that divide the numerator of $\nu(u)$
for some $u\in U$ with $\psi(u)\ne 1$.

This gives a contradiction, completing the proof.
\end{pf}

\section{Double cosets}\label{doublecosets}
In this section we use Theorem \ref{theoremA} to prove Theorem \ref{theoremC}.

The group $G$ in Theorem \ref{theoremC} is virtually soluble.
Indeed, one may immediately reduce to the case where $G$ is
soluble, since the hypotheses and the conclusion of 
Theorem \ref{theoremC} are stable under the passage to finite
index subgroups.

Similarly, we may replace the finitely generated virtually abelian
subgroups $A,B$ of Theorem \ref{theoremC} by any finite-index subgroups
of $A,B$ respectively.  In particular we may assume that $A,B$ are free
abelian of finite rank.

We argue by double induction, on the derived length of the soluble
group $G$ and on the sum of the ranks of the free abelian groups
$A,B$.

If $G$ is abelian, then the hypotheses imply that $G$ is finitely
generated (for example, by the union of a basis for $A$, a basis for $B$,
and a set of double-coset representatives).  Otherwise,
$G$ has an abelian normal subgroup $N$ such that the derived length of $G/N$
is less than that of $G$.  Since $|AN\backslash G/BN|<\infty$, it
follows by inductive hypothesis that $G/N$ is virtually polycyclic,
and it suffices to show that $N$ is finitely generated as an abelian
group.

Now $|A\backslash AN/(B\cap AN)|<\infty$, so there is no loss of generality
in supposing that $G=AN$.  Then $A\cap N$ is central, and hence normal
in $G$.  If $A\cap N\ne\{1\}$, then we may apply the inductive hypothesis
together with the fact that
$$\left|\frac{A}{A\cap N} \left\backslash\frac{G}{A\cap N}\right/\frac{B\cdot(A\cap N)}{A\cap N}\right|<\infty$$
to deduce that $G/(A\cap N)$ is virtually polycyclic, and hence that $G$ is virtually
polycyclic.  This reduces us to the case where $A\cap N=\{1\}$.

In a similar way, we may assume that $BN=G$ and that $B\cap N=\{1\}$.  Hence
$G=N\rtimes A=N\rtimes B$.  In particular, $A\cong B$
via an isomorphism $\theta:A\to B$ uniquely deterined by the property that
$a\theta(a)^{-1}\in N$ for all $a\in A$.  Moreover, $N$ is a $\Z A$-module,
and the map $\delta:A\to N$ defined by $\delta(a)=a\theta(a)^{-1}$ is
a derivation.  For $g\in G$, the double coset $AgB$ intersects $N$: indeed
if $g\in N$ then $AgB\cap N$ is precisely the orbit of $A$ under the 
affine action $$a*g=a\cdot g+\delta(a)=aga^{-1}.a\theta(a)^{-1}=ag\theta(a)^{-1}.$$

In particular, $N$ is the union of finitely many $A$-orbits under this
action, so $N$ is finitely generated as an abelian group, by Theorem \ref{theoremA}.
This concludes the proof of Theorem \ref{theoremC}.

\end{document}